\documentclass[%
 aip,
 amsmath,amssymb,
 reprint,%
]{revtex4-1}

\usepackage{graphicx}
\usepackage{dcolumn}
\usepackage{bm}

\usepackage[utf8]{inputenc}
\usepackage[T1]{fontenc}
\usepackage{mathptmx}

\usepackage{subcaption} 
\usepackage{nicefrac} 
\usepackage{multirow} 

\usepackage{xcolor}
\usepackage{ulem} 

\begin{document}

\preprint{AIP/123-QED}

\title[Lorenz System Stability Identification]{Lorenz System State Stability Identification using Neural Networks}

\author{Megha Subramanian}
 \email{megha.subramanian@pnnl.gov}
\author{Ramakrishna Tipireddy}%
 \email{ramakrishna.tipireddy@pnnl.gov.}
\author{Samrat Chatterjee}%
 \email{samrat.chatterjee@pnnl.gov.}
\affiliation{ 
Pacific Northwest National Laboratory, Richland, WA
}%


\date{\today}

\begin{abstract}
Nonlinear dynamical systems such as Lorenz63 equations are known to be chaotic in nature and sensitive to initial conditions. As a result, a small perturbation in the initial conditions results in deviation in state trajectory after a few time steps. System identification is often challenging in such systems especially when the solution transitions from one regime to another compared to identification when the trajectory of the solution resides within the same regime. The algorithms and computational resources needed to accurately identify the system states vary depending on whether the solution is in transition region or not. We refer to the transition and non-transition regions as unstable and stable regions respectively. We label a system state to be stable if it's immediate past and future states reside in the same regime. However, at a given time step we don't have the prior knowledge about whether system is in stable or unstable region. In this paper, we develop and train a feed forward (multi-layer perceptron) Neural Network to classify the system states of a Lorenz system as stable and unstable. We pose this task as a supervised learning problem where we train the neural network on Lorenz system which have states labeled as stable or unstable. We then test the ability of the neural network models to identify the stable and unstable states on a different Lorenz system that is generated using different initial conditions. We also evaluate the classification performance in the mismatched case i.e., when the initial conditions for training and validation data are sampled from different intervals. We show that certain normalization schemes can greatly improve the performance of neural networks in especially these mismatched scenarios. The classification framework developed in the paper can be a preprocessor for a larger context of sequential decision making framework where the decision making is performed based on observed stable or unstable states.
\end{abstract}

\maketitle

\section{\label{sec:Introduction} Introduction}






Automated decision making in complex systems relies on robust system state identification under uncertainty. System state identification via data-driven learning is the often the first step within a sequential decision making (SDM) framework \cite{amiri2019robot}. An SDM framework may comprise of multiple computational modules for system state estimation, context-aware reasoning, and probabilistic planning. A computational agent in such settings may continuously interact with an external environment, compute the best policies (state to action mappings), and execute actions that maximize learning reward signals in the long run \cite{sutton2018reinforcement}. In this context, system state is often hidden or partially observable from an agent. This requires an agent to infer the state of the system indirectly through a combination of data driven learning and domain knowledge. Typically the system state space can be characterized in discrete terms to represent different operating regimes. A discrete state space representation may result in translation of the state estimation problem as a data-driven classification task. In recent years, use of neural network classifiers have led to high accuracy and performance across multi-modal data streams (for example text and image)\cite{he2015deep}. 
State estimation of complex systems such as nonlinear dynamical systems (for example Lorenz63 equations) is challenging due to fluctuating system trajectories and high sensitivity to initial conditions. In this paper, we demonstrate a proof of concept application with a neural network classifier for state estimation of a nonlinear chaotic dynamical system. 

Neural networks are a sequence of densely interconnected nodes. Through a process known as training, the networks learn the parameters, or the relationship between these nodes. The availability of large amounts of data and massive computational power are some of the reasons that have facilitated training of deep neural networks in recent years. Consequently, deep neural networks have garnered widespread attention across different scientific disciplines including speech recognition\cite{6296526}, image classification\cite{he2015deep}, language modeling\cite{NIPS2017_3f5ee243} and even protein folding\cite{alma9926190102476}. The idea of leveraging neural networks to study dynamical systems has been explored for many years\cite{alma9927832002476}. Recently, both feed-forward and recurrent architectures have been used for several forecasting and prediction based tasks on chaotic systems. For example, feed-forward neural networks have been shown to predict extreme events in H\'enon map \cite{Lellep_2020} and LSTM\cite{hochreiter1997long} architectures have been used for forecasting high-dimensional chaotic systems \cite{alma9929220102476}. Additionally Reservoir Computing based approaches have been leveraged for data-driven prediction of chaotic systems (Refs.~\onlinecite{article,DBLP:journals/corr/abs-1803-04779,PhysRevLett.120.024102}). There is also literature suggesting that a hybrid approach to forecasting involving both machine learning and knowledge based models \cite{DBLP:journals/corr/abs-1803-04779} leads to higher prediction accuracy for a longer period on chaotic systems. Another line of work involves discovering underlying models from data, wherein autoencoder network architectures are used to recognize the coordinate transformation and the governing equations of the dynamical systems \cite{champion2019data}. All these examples highlight the general feasibility of using machine learning and neural networks to study a wide variety of tasks involving chaotic systems. 
Previous work  has shown that feed-forward neural networks are reasonable candidates for predicting regime changes and duration in Lorenz63 systems \cite{brugnago2020classification}. Inspired by this approach, we leverage feed-forward neural networks for our novel task of classifying the data points of Lorenz system. The solution of the Lorenz system consists of two regimes which we refer to as the left and right regimes. Unlike the approach proposed in \onlinecite{brugnago2020classification}, we do not restrict ourselves to the boundary $x=0$ to distinguish between the left and right regimes. Since we are interested in studying a broader range of Lorenz systems, we adopt the mean value of x-coordinates as a boundary to distinguish between the left and right regimes. Additionally, unlike \onlinecite{brugnago2020classification}, we do not constrain ourselves to the interval $\left(0,1\right)$ for sampling the initial conditions for the variables $\left(x,y,z\right)$ for training and validation datasets. Instead, we sample the initial conditions for training and validation datasets from a wide range of intervals (Sections ~\ref{subsec : Initial Conditions for Training Dataset} and ~\ref{subsec:Initial Conditions for Validation Dataset}) and analyze the behaviour of the classifier in these intervals.

One common problem associated with neural network architectures is lack of generalization with validation data that is drawn from a different distribution as compared to the data used for training. This mismatch between training and validation distributions, also referred to as covariate shift \cite{sugiyama2012machine} is a well-studied aspect in most of the scientific fields that have benefited from advances in neural network architectures. For example, within the speech recognition community, this mismatch is referred to as the acoustic mismatch and is one of the most critical factors that affects the deployment of speech recognition systems in real world environments \cite{VINCENT2017535}. A similar terminology used in the field of natural language processing is domain mismatch, which renders tasks like cross-lingual document classification especially challenging \cite{lai2019bridging}. In the context of chaotic systems, there has been some published work \cite{scher2019generalization} that discusses the generalization capabilities of neural networks trained on Lorenz systems. One of the results of this work on Lorenz63 systems is that neural networks trained only on a part of the system's phase space struggle to skillfully forecast in regions that are excluded from the training phase space. These observations are consistent with our initial results in Section \ref{subsec : Initial Conditions for Training and Validation Data sampled from different intervals}, where we notice and analyze the drop in the performance of neural networks in the mismatched case.

In this paper, we further show how a very simple yet effective normalization approach helps in counteracting the mismatch problem. As we will show in Section \ref{subsec:Normalization}, this pre-processing step of normalization helps in training neural networks that give an impressive performance on a wide variety of validation data for the task of classification of stable/unstable data points. To the best of our knowledge, leveraging neural networks to classify data points of Lorenz system as stable/unstable, analyzing their performance in the mismatched case and the proposed normalization scheme have not been published so far.

The rest of the paper is organized as follows : Section \ref{sec:Lorenz System of Equations} gives an overview of the Lorenz system of equations used in our study. The Section \ref{sec:Labeling Strategy} details the labeling strategy used for our supervised learning problem. The Section \ref{sec:Training and Validation Datasets} explains how the training and validation datasets are generated. The architecture of neural network is described in Section \ref{sec:Neural Network Architecture}. The analysis of initial results without the normalization are presented in Section \ref{sec:Results}. The Section \ref{subsec:Normalization} gives more insight into the normalization procedure and the associated results. The Section \ref{sec:Conclusion} discusses the future scope of work.

\section{\label{sec:Lorenz System of Equations} Lorenz System of Equations}
Lorenz System is a dynamical system consisting of three coupled differential equations. 
\begin{subequations}
\begin{equation}
\label{subeqn:Lorenz X Derivative}
    \frac{dx}{dt} = \sigma(y-x)
\end{equation}
\begin{equation}
\label{subeqn:Lorenz Y Derivative}
    \frac{dy}{dt} = x(\rho - z) - y
\end{equation}
\begin{equation}
\label{subeqn:Lorenz Z Derivative}
    \frac{dz}{dt} = xy - \beta z
\end{equation}
\end{subequations}
where $x, y, z$ are the state variables and $\sigma$, $\rho$, $\beta$ are the parameters. In this work, we set $\sigma = 10$, $\rho = 28$ and $\beta = \nicefrac{8}{3}$ \cite{lorenz1963deterministic}. We represent the time derivatives by $\frac{dx}{dt}$, $\frac{dy}{dt}$ and $\frac{dz}{dt}$  respectively and $x_{0}$, $y_{0}$ and $z_{0}$ are the initial conditions of the Lorenz system in $x$, $y$ and $z$ directions respectively.

\begin{figure}\label{fig:lorenz}
\includegraphics[width=0.5\textwidth]{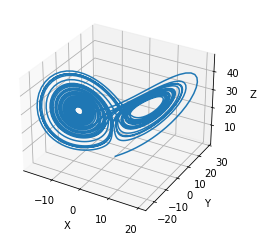}
\caption{A typical Lorenz System}
\label{fig:A typical Lorenz System}
\end{figure}

Fig.~\ref{fig:A typical Lorenz System} shows the trajectory of system states of a typical Lorenz system. The Lorenz system is a highly chaotic system. Even a small perturbation in the initial conditions can cause the trajectories to diverge significantly. The system is characterized by two distinctive regimes that give rise to a butterfly shaped figure. Starting from some initial conditions, the data points of a Lorenz system can either remain in the same regime or can alternate between the two regimes. In this work, we are interested in identifying those data points of the Lorenz system that will either undergo a regime change in the future time steps or have undergone a regime change in the past few time steps. In other words, we are interested in isolating those data points whose past and future data points lie in different regimes. We refer to such data points as unstable data points. We train a  neural network to distinguish between stable and unstable data points based on labeled training examples of Lorenz systems. We then test the ability of the neural network models to identify the unstable data points on unseen validation data. Since the Lorenz system is a highly chaotic system, even a small change in the initial conditions will cause a major shift in the location of stable and unstable data points. Consequently, it is a challenging task for neural network models to isolate the unstable data points on unseen validation data.

In this work, we use neural network models to classify unstable data points in both matched (i.e when the initial conditions for training and validation data are sampled from the same interval) and mismatched conditions (i.e when the initial conditions for training and validation data are sampled from different intervals). We first demonstrate that it is particularly difficult for neural network models to reliably classify the unstable data points in mismatched conditions. We further observed that certain normalization schemes can greatly improve the performance of neural network models in mismatched conditions. 
\section{\label{sec:Labeling Strategy} Labeling Strategy}
The classification of data points as stable or unstable can be posed a supervised learning problem that requires labeled data points during training. We use a two-step process to generate the labels for stable/unstable classification. In the first step, we introduce a heuristic method to identify the right and left regimes of the Lorenz attractor. In the second step, we label the data points as stable or unstable based on the regimes identified in the first step. 
\subsection{\label{subsec:Identification of Left and Right Regimes}Identification of Left and Right Regimes}
\begin{figure}
\includegraphics[width=0.5\textwidth]{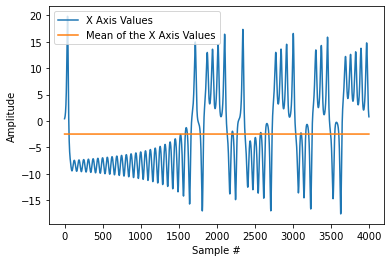}
\caption{\label{fig:X Axis Values and their Mean}X Axis Values and their Mean}
\end{figure}

\begin{figure}
\includegraphics[width=0.5\textwidth]{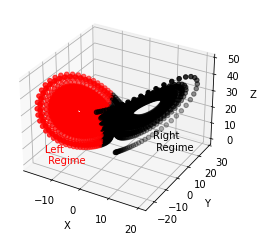}
\caption{Lorenz System with labelled Left and Right Regime.The regimes are labelled based on mean value of the x-coordinate}
\label{fig:Lorenz System with Left and Right Regime}
\end{figure}

In \onlinecite{brugnago2020classification}, the authors sampled the initial
conditions for the dynamic variables from the interval $\left[0, 1\right]$ for their experiments. They used the condition $x < 0$ to define the left regimes and $ x \geq 0$ to define the right regimes. Since we want to apply the results of our classification to a wide variety of Lorenz systems which may not be perfectly aligned with the boundary $x=0$, we adopt a different approach. Specifically, we define the left and right regimes of the Lorenz attractor based on the mean value of x coordinate. The Fig.~\ref{fig:X Axis Values and their Mean} shows the x-axis values and their mean.\\
In our heuristic-approach, if the x coordinate of the data point has a value greater than the mean in the x direction, ($ x > mean(x) $), the data point belongs to the right regime, otherwise it belongs to the left regime ($ x < mean(x)$). The Fig.~\ref{fig:Lorenz System with Left and Right Regime} shows a Lorenz System with left and right regimes labeled using this approach.\\
\subsection{\label{subsec:Labeling the data points as stable or unstable}Labeling the data points as stable or unstable}
\begin{figure}
\includegraphics[width=0.5\textwidth]{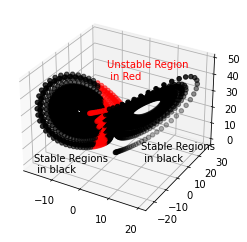}
\caption{Lorenz System with Stable and Unstable Regions}
\label{fig:Lorenz System with Stable and Unstable Regions}
\end{figure}

After defining the left and right regimes, we label the data points at every time step $t$ as stable or unstable. Specifically, we consider a window around the data point at time step $t$. In this study, we work with a window of $5$ data points from the past time steps $ t-1,\; t-2,\; t-3,\; t-4,\; t-5$ and $5$ data points from the future time steps $t+1,\; t+2,\; t+3,\; t+4,\; t+5$. A data point at time step $t$ is considered stable if the past $5$ and the future $5$ neighbouring data points belong to the same regime. The data point is labeled unstable even if one of the $10$ neighbours belong to the different regime. Here the number of points is a hyper parameter that we chose as based on trial and error to manually classify the regions into stable and unstable. This number can be calibrated in advance for different dynamical systems. The Fig.~\ref{fig:Lorenz System with Stable and Unstable Regions} shows a Lorenz system with labeled stable and unstable data points.

\section{\label{sec:Training and Validation Datasets} Training and Validation Datasets}
For generating a Lorenz System such as the one shown in Fig.~\ref{fig:A typical Lorenz System}, we specify the initial condition 
 $x_{0}$, $y_{0}$ and $z_{0}$ for the dynamic variables. We then integrate the system with a time step of $0.01$ using the ODEINT\cite{Ahnert_2011} with Runge-Kutta 45 (RK45) algorithm provided by scipy\cite{2020SciPy-NMeth}. For every Lorenz System such as the one shown in Fig.~\ref{fig:Lorenz System with Stable and Unstable Regions}, we generate $4000$ data points. For training, we use $25$ such Lorenz Systems, resulting in a total of $100,000$ data points for the training set. For validation, we use $5$ such Lorenz Systems, resulting in a total of $20,000$ data points for the validation set.
\subsection{\label{subsec : Initial Conditions for Training Dataset}Initial Conditions for Training Dataset}
We experiment with three different initial condition intervals for the training data set. The Eqs.~(\ref{eq:Initial Conditions Train}) describe the intervals from which the initial condition values for the training data are sampled randomly. It is important to remember that each of the initial conditions $x_{0}$, $y_{0}$ and $z_{0}$ are sampled independently from the intervals. 
\begin{subequations}
\label{eq:Initial Conditions Train}
\begin{equation}
\label{subeq : Initial Conditions Train 1}
    x_{0}, y_{0}, z_{0}\sim \left[0, 1\right]
\end{equation}
\begin{equation}
\label{subeq : Initial Conditions Train 2}
    x_{0}, y_{0}, z_{0}\sim \left[-1, 0\right]
\end{equation}
\begin{equation}
\label{subeq : Initial Conditions Train 3}
    x_{0}, y_{0}, z_{0}\sim \left[-1, 1\right]
\end{equation}
\end{subequations}
For each of the initial conditions shown in Eqs.~(\ref{eq:Initial Conditions Train}), we generate $25$ Lorenz Systems according to the procedure described above. These Lorenz systems are then used to train different neural networks. We are interested in evaluating whether a neural network trained on chaotic Lorenz system with the initial condition intervals specified in Eqs.~(\ref{eq:Initial Conditions Train}) can achieve accurate classification performance on validation sets generated from a wide variety of initial conditions.
\subsection{\label{subsec:Initial Conditions for Validation Dataset}Initial Conditions for Validation Dataset}
The Eqs.~(\ref{eq:Initial Conditions Val}) describe the intervals from which the initial condition values for the validation data are sampled randomly for our experiments. We use different random initial seeds while generating the initial conditions for training and validation data. This ensures that the randomly sampled initial conditions for training and validation data are different, even though the intervals from which these are sampled might be the same. 
\begin{subequations}
\label{eq:Initial Conditions Val}
\begin{equation}
\label{subeq : Initial Conditions Val 1}
    x_{0}, y_{0}, z_{0}\sim \left[0, 1\right]
\end{equation}
\begin{equation}
\label{subeq : Initial Conditions Val 2}
    x_{0}, y_{0}, z_{0}\sim \left[-1, 0\right]
\end{equation}
\begin{equation}
\label{subeq : Initial Conditions Val 3}
    x_{0}, y_{0}, z_{0}\sim \left[-1, 1\right]
\end{equation}
\begin{equation}
\label{subeq : Initial Conditions Val 4}
    x_{0}, y_{0}, z_{0}\sim \left[2, 4\right]
\end{equation}
\begin{equation}
\label{subeq : Initial Conditions Val 5}
    x_{0}, y_{0}, z_{0}\sim \left[0, 10\right]
\end{equation}
\begin{equation}
\label{subeq : Initial Conditions Val 6}
    x_{0}, y_{0}, z_{0}\sim \left[-10, 10\right]
\end{equation}
\end{subequations}

As we will see in section ~\ref{sec:Results}, it is particularly challenging for neural networks to accurately classify the stable and unstable points when the initial conditions for training and validation data are sampled from different intervals. We will also see how the normalization scheme discussed in section ~\ref{subsec:Normalization} can help improve the performance in these cases.

\section{\label{sec:Neural Network Architecture} Neural Network Architecture}
\begin{figure}
\includegraphics[width=0.5\textwidth]{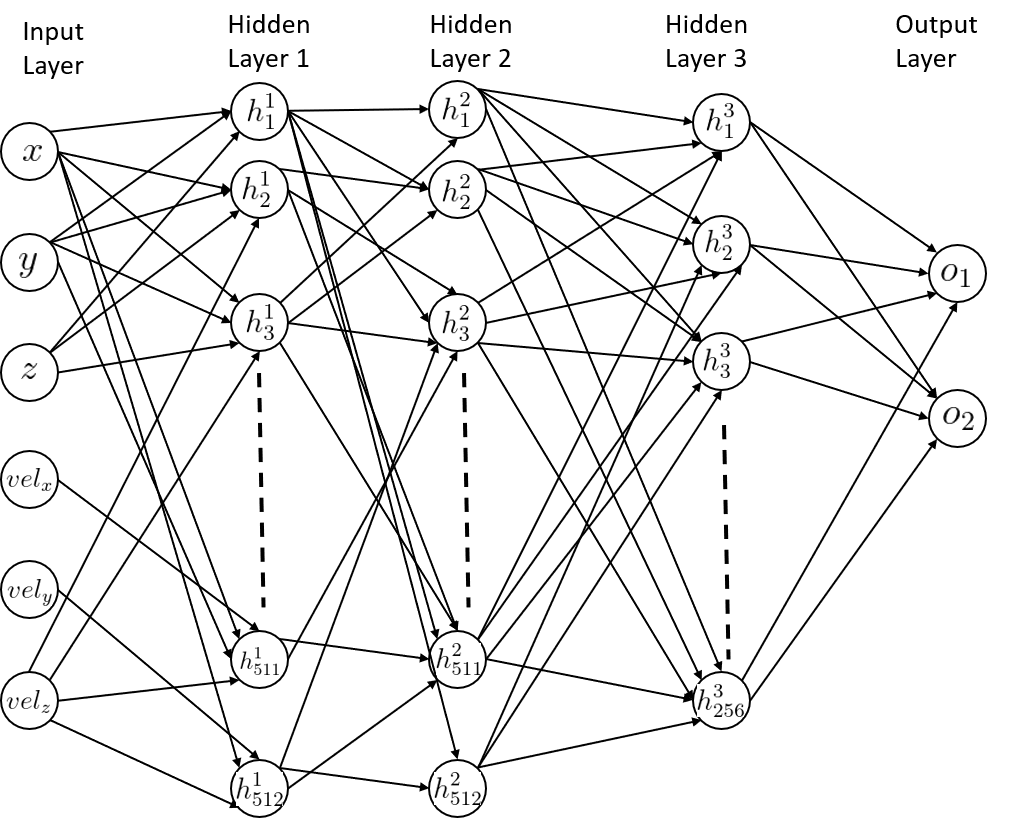}
\caption{Fully Connected Feed Forward Neural Network Architecture used in the experiments. Note that $h_{m}^{n}$ denotes $m^{th}$ hidden node in $n^{th}$ layer}
\label{fig:Neural Network Architecture}
\end{figure}

In this section, we describe the neural network architecture used to classify the data points of the Lorenz System (See Fig. \ref{fig:Neural Network Architecture}). The neural network takes as input the three Lorenz variables  $ x, y, z$ of the data point as well as the three time derivatives $\frac{dx}{dt}$, $\frac{dy}{dt}$, $\frac{dz}{dt}$. Thus, the input to the network is a $6$ dimensional feature vector. The neural network outputs the probabilities of the input feature vector being stable or unstable. The network consists of four fully connected layers. The first two layers have $512$ neurons each and have hyperbolic tangents and rectified linear unit (ReLU) activation functions respectively. The third layer consists of $256$ neurons and ReLU activation function. The final layer consists of $2$ neurons and sigmoid activation function that predicts the probabilities of the data points belonging to
either the stable or unstable classes. We use Adam optimizer \cite{kingma2017adam} with a learning rate of $0.001$ and binary cross entropy loss function. 
The neural network was implemented using the Keras library\cite{chollet2015keras}. 
\section{\label{sec:Results} Results}
As mentioned above, the Lorenz system is a chaotic system which is highly sensitive to initial conditions. Even a small perturbation in the initial conditions can cause the trajectories to diverge to a large extent and can affect the location of stable and unstable data points. Hence, it is a non-trivial task for neural network models to predict the unstable data points of Lorenz Systems which have not been encountered during training. As we discuss in succeeding sections, it is especially challenging for neural networks to predict the unstable data points of Lorenz systems in mismatched conditions i.e when the initial conditions of validation Lorenz systems have been sampled from a different interval as compared to Lorenz systems used for training. In subsection~\ref{subsec : Initial Conditions for Training and Validation Data sampled from same intervals}, we present the results of predicting unstable data points using neural networks when the initial conditions for training and validation have been sampled from the same interval. In subsection~\ref{subsec : Initial Conditions for Training and Validation Data sampled from different intervals}, we discuss the results in mismatched conditions. Note that our dataset is imbalanced. The number of stable data points (majority class) are much more than unstable data points (minority class). In such cases, using the accuracy scores as performance evaluation metric is misleading since the model will have a high accuracy even if it always predicts the majority class. Hence we do not report our results in terms of accuracy score. Instead we report the results in terms of precision and recall scores which are defined in Eqs.~\ref{subeq:Precision and Recall Definitions}.

\begin{subequations}
\label{subeq:Precision and Recall Definitions}
\begin{equation}
    \label{subeq : Precision Definition}
    Precision = \frac{True Positives}{True Positives + False Positives}
\end{equation}
\begin{equation}
    \label{subeq : Recall Definition}
    Recall = \frac{True Positives}{True Positives + False Negatives}
\end{equation}

\end{subequations}

\subsection{\label{subsec : Initial Conditions for Training and Validation Data sampled from same intervals}Initial Conditions for Training \& Validation Data sampled from same interval}
\begin{figure}[!tbp]
    \begin{subfigure}[b]{0.49\columnwidth}
        \includegraphics[width = 0.99\columnwidth]{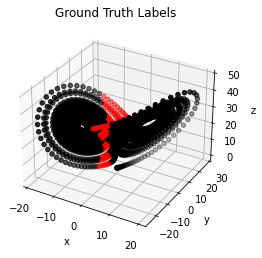}
        \caption{Ground Truth Labels}
        \label{fig:LS_Res1_Predicted}
    \end{subfigure}
    \begin{subfigure}[b]{0.49\columnwidth}
        \includegraphics[width = 0.99\columnwidth]{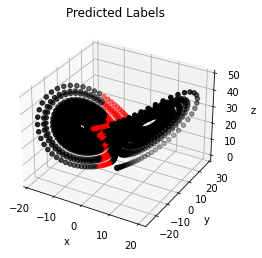}
        \caption{Predicted Labels}
        \label{fig:LS_Res1_GT}
    \end{subfigure}
\caption{Classification result using neural network models in matched conditions. The initial conditions for training data were sampled according to Eq.~(\ref{subeq : Initial Conditions Train 1}) and for validation data were sampled according to Eq.~(\ref{subeq : Initial Conditions Val 1}). The precision and recall scores are $0.9$ and $0.744$ respectively}
\label{fig : Res1}
\end{figure}

\begin{figure}[!tbp]
    \begin{subfigure}[b]{0.49\columnwidth}
        \includegraphics[width = 0.99\columnwidth]{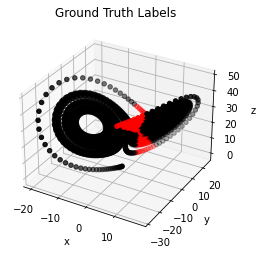}
        \caption{Ground Truth Labels}
        \label{fig:LS_Res2_Predicted}
    \end{subfigure}
    \begin{subfigure}[b]{0.49\columnwidth}
        \includegraphics[width = 0.99\columnwidth]{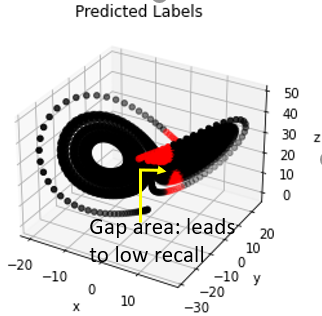}
        \caption{Predicted Labels}
        \label{fig:LS_Res2_GT}
    \end{subfigure}
\caption{Classification result using neural network models in matched conditions. The initial conditions for training data were sampled according to Eq.~(\ref{subeq : Initial Conditions Train 2}) and for validation data were sampled according to Eq.~(\ref{subeq : Initial Conditions Val 2}). The precision and recall scores are $0.94$ and $0.78$ respectively}
\label{fig : Res2}
\end{figure}

The example shown in the Fig~\ref{fig : Res1} shows the results of applying neural networks for predicting unstable data points when the initial condition variables $x_{0}$, $y_{0}$ and $z_{0}$ for both the training and the validation data are sampled from the same interval. The initial conditions for the training data for the model are sampled according to Eq.~(\ref{subeq : Initial Conditions Train 1}), and those for validation data are sampled according to Eq.~(\ref{subeq : Initial Conditions Val 1}). It is important to note that even though the sampling intervals for initial conditions are the same, the randomly sampled initial conditions for training and validation data within this interval are different. We observe that the neural network performs reasonably well in predicting the unstable data points with a precision of $0.90$ and a recall of $0.744$. The Fig~\ref{fig : Res2} shows another example of applying neural networks for classification of unstable data points when the initial conditions for both the training and validation data are sampled from the same interval. In this example, the initial conditions for the training data for the model are sampled according to Eq.~(\ref{subeq : Initial Conditions Train 2}) and those for validation data are sampled according to Eq.~(\ref{subeq : Initial Conditions Val 2}). The precision and recall scores obtained are $0.94$ and $0.78$ respectively.\\
We observe that for both the examples, the recall scores are lower than those of precision. According to Eq.~(\ref{subeq : Recall Definition}), a lower recall score corresponds to more number of false negatives. This effect is better seen in Fig~\ref{fig : Res2}, where the model sometimes fails to predict the unstable data points in the gap region. This observation is also prevalent in mismatched conditions. We will see in section~\ref{subsec:Normalization} that the normalization scheme helps in reducing the false negatives and in improving the recall scores. 

\subsection{\label{subsec : Initial Conditions for Training and Validation Data sampled from different intervals}Initial Conditions for Training \& Validation Data sampled from different intervals}
\begin{figure}[!tbp]
 \begin{subfigure}[b]{0.49\columnwidth}
        \includegraphics[width = 0.99\columnwidth]{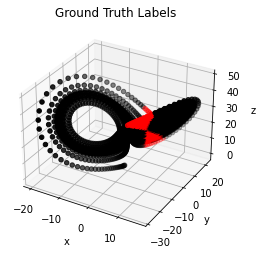}
        \caption{Ground Truth Labels}
        \label{fig:LS_Res3_GT}
    \end{subfigure}
    \begin{subfigure}[b]{0.49\columnwidth}
        \includegraphics[width = 0.99\columnwidth]{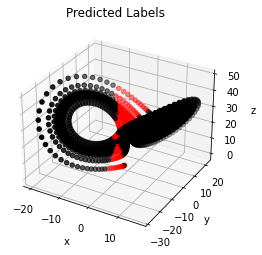}
        \caption{Predicted Labels}
        \label{fig:LS_Res3_Predicted}
    \end{subfigure}
    \caption{Classification result using neural network models in mismatched conditions. The initial conditions for training data were sampled according to Eq.~(\ref{subeq : Initial Conditions Train 1}) and for validation data were sampled according to Eq.~(\ref{subeq : Initial Conditions Val 2}). The precision and recall scores are $0.01$ and $0.004$ respectively}
    \label{fig:Res3}
\end{figure}

\begin{figure}[!tbp]
 \begin{subfigure}[b]{0.49\columnwidth}
        \includegraphics[width = 0.99\columnwidth]{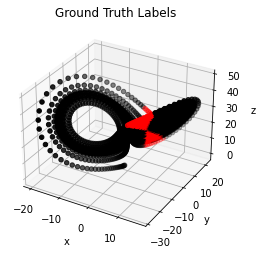}
        \caption{Ground Truth Labels}
        \label{fig:LS_Res4_GT}
    \end{subfigure}
    \begin{subfigure}[b]{0.49\columnwidth}
        \includegraphics[width = 0.99\columnwidth]{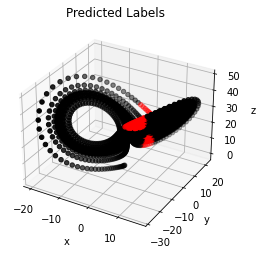}
        \caption{Predicted Labels}
        \label{fig:LS_Res4_Predicted}
    \end{subfigure}
    \caption{Classification result using neural network models in matched conditions. The initial conditions for training data were sampled according to Eq.~(\ref{subeq : Initial Conditions Train 2}) and for validation data were sampled according to Eq.~(\ref{subeq : Initial Conditions Val 2}). The precision and recall scores are $0.7955$ and $0.5659$ respectively.}
    \label{fig:Res4}
\end{figure}

\begin{figure}[!tbp]
 \begin{subfigure}[b]{0.49\columnwidth}
        \includegraphics[width = 0.99\columnwidth]{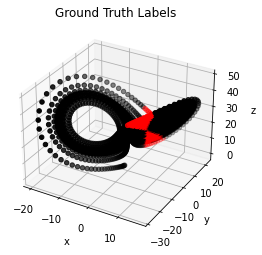}
        \caption{Ground Truth Labels}
        \label{fig:LS_Res5_GT}
    \end{subfigure}
    \begin{subfigure}[b]{0.49\columnwidth}
        \includegraphics[width = 0.99\columnwidth]{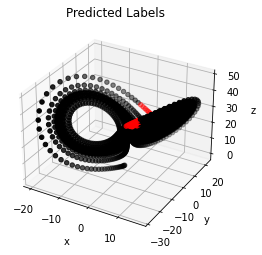}
        \caption{Predicted Labels}
        \label{fig:LS_Res5_Predicted}
    \end{subfigure}
    \caption{Classification result using neural network models in mismatched conditions. The initial conditions for training data were sampled according to Eq.~(\ref{subeq : Initial Conditions Train 3}) and for validation data were sampled according to Eq.~(\ref{subeq : Initial Conditions Val 2}). The precision and recall scores are $0.846$ and $0.125$ respectively}
    \label{fig:Res5}
\end{figure}

The Fig.~\ref{fig:Res3} illustrates the performance of neural network models in mismatched condition. The initial conditions for training data are sampled according to Eq.~(\ref{subeq : Initial Conditions Train 1}) and for validation data are sampled according to Eq.~(\ref{subeq : Initial Conditions Val 2}). In this example, the neural network models almost completely miss the region of unstable data points and fail to identify them reliably. This is also reflected in the low precision and recall values of $0.01$ and $0.004$ respectively. As a reference, the Fig.~\ref{fig:Res4} shows the result on the same validation set as Fig.~\ref{fig:Res3} using a neural network model that has the initial conditions for training data sampled according to Eq.~\ref{subeq : Initial Conditions Train 2}. A comparison of Fig.~\ref{fig:Res3} and Fig.~\ref{fig:Res4}, shows the importance of initial conditions on the performance of neural networks for the specific task of identification of unstable data points of Lorenz system. One can clearly see that using the features described in Section~\ref{sec:Neural Network Architecture} as-is, limits the usability of the neural network models, in that, the models cannot be used for the mismatched scenario.\\ 
Neural network models that are generalizable should be able to perform reliably on a broad variety of validation data. To this end, we train neural network models that have the initial conditions for the training data sampled from a slightly larger interval. Specifically, we use Eq.~(\ref{subeq : Initial Conditions Train 3}) to sample the initial conditions for training data. The Fig.~\ref{fig:Res5} shows the classification result when the initial conditions for training data were sampled according to Eq.~(\ref{subeq : Initial Conditions Train 3}) and for validation data were sampled according to Eq.~(\ref{subeq : Initial Conditions Val 2}). The precision and recall values are $0.846$ and $0.125$ respectively. Again, the low recall score indicates a high number of false negatives (which is also reflected as less number of data points labelled in red in Fig.~\ref{fig:Res5}). Although the classification results of the models shown in Fig.~\ref{fig:Res5} are not comparable to matched case shown in Fig.~\ref{fig:Res4}, these models are able to capture the region of unstable data points better than those in Fig.~\ref{fig:Res3}. Since it is impractical to have a matched model for every validation set (one whose initial conditions during training have been sampled from the same interval, such as the one shown in Fig.~\ref{fig:Res4}), we think that model trained according to Eq.~(\ref{subeq : Initial Conditions Train 3}) (such as the one shown in Fig.~\ref{fig:Res5}) might serve as a good candidate for training neural networks that can perform reliably on a wide variety of validation data. In the section ~\ref{subsec:Normalization}, we show how the normalization scheme helps in further improving the performance of models trained according to Eq.~(\ref{subeq : Initial Conditions Train 3}) in mismatched conditions.



\section{\label{subsec:Normalization}Normalization}
\begin{figure}[!tbp]
    \begin{subfigure}[b]{0.49\columnwidth}
        \includegraphics[width = 0.99\columnwidth]{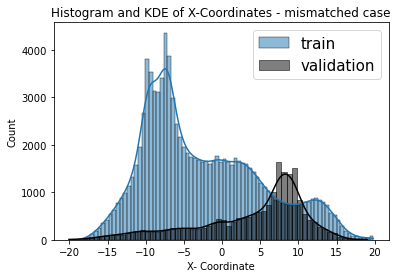}
        \caption{Histogram \& KDE - mismatched case}
        \label{fig:Histogram of X-Coordinates - mismatched case}
    \end{subfigure}
    \begin{subfigure}[b]{0.49\columnwidth}
        \includegraphics[width = 0.99\columnwidth]{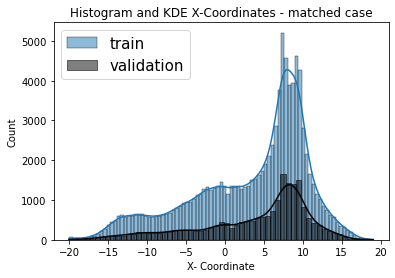}
        \caption{Histogram \& KDE - matched case}
        \label{fig : Histogram of X-Coordinates - matched case}
    \end{subfigure}
\caption{Histogram and Kernel Density Estimates 
(KDE) of X-Coordinates of the training and validation data for the mismatched and matched cases. The initial conditions for training data were sampled according to Eq.~(\ref{subeq : Initial Conditions Train 1}) for Fig.~\ref{fig:Histogram of X-Coordinates - mismatched case} and according to Eq.~(\ref{subeq : Initial Conditions Train 2}) for Fig.~\ref{fig : Histogram of X-Coordinates - matched case}. The initial conditions for validation data were sampled according to Eq.~(\ref{subeq : Initial Conditions Val 2}) for both the cases}
\label{Fig : Histogram Comparisons - Before Normalization}
\end{figure}

\begin{figure}[!tbp]
    \begin{subfigure}[b]{0.49\columnwidth}
        \includegraphics[width = 0.99\columnwidth]{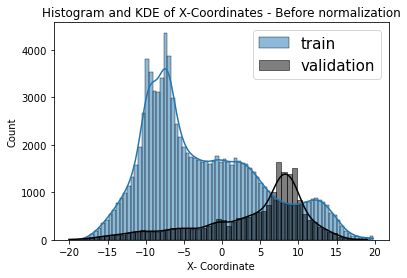}
        \caption{Histogram \& KDE - mismatched case before normalization}
        \label{fig:Histogram of X-Coordinates - mismatched case before normalization}
    \end{subfigure}
    \begin{subfigure}[b]{0.49\columnwidth}
        \includegraphics[width = 0.99\columnwidth]{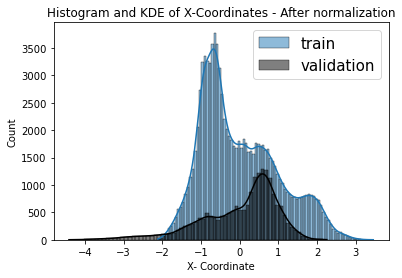}
        \caption{Histogram \& KDE - mismatched case after normalization}
        \label{fig : Histogram of X-Coordinates - mismatched case after normalization}
    \end{subfigure}
\caption{Histogram and Kernel Density Estimates (KDE) of X-Coordinates of the training and validation data for the mismatched case before and after normalization. The initial conditions for training data for Fig.~\ref{fig:Histogram of X-Coordinates - mismatched case before normalization} and Fig.~\ref{fig : Histogram of X-Coordinates - mismatched case after normalization} were sampled according to Eq.~(\ref{subeq : Initial Conditions Train 1}). The initial conditions for validation data were sampled according to Eq.~(\ref{subeq : Initial Conditions Val 2})} 
\label{Fig : Histogram Comparisons - Before and After Normalization}
\end{figure}

In order to understand why the classification results shown in Fig.~\ref{fig:Res3} are worse than those in Fig.~\ref{fig:Res4}, we look at the statistics of the training and validation data. The Fig.~\ref{Fig : Histogram Comparisons - Before Normalization} compares the histograms and kernel density estimates of the features of training and validation data in the mismatched and matched cases. For simplicity, only the x-coordinate features are shown. In the Fig.~\ref{fig:Histogram of X-Coordinates - mismatched case}, the initial conditions of training data were sampled according to Eq.~(\ref{subeq : Initial Conditions Train 1}) and those of the validation data were sampled according to Eq.~(\ref{subeq : Initial Conditions Val 2}). One can notice that there is a strong mismatch between the histograms and kernel density estimates of training and validation data in Fig.~\ref{fig:Histogram of X-Coordinates - mismatched case}. This mismatch leads to a degradation in the performance of the neural networks. In contrast, the initial conditions of the training data were sampled according to Eq.~(\ref{subeq : Initial Conditions Train 2}) in Fig.~\ref{fig : Histogram of X-Coordinates - matched case}. The initial conditions for validation data were sampled according to Eq.~(\ref{subeq : Initial Conditions Val 2}). Since the initial conditions for training and validation data were sampled from the same interval in Fig.~\ref{fig : Histogram of X-Coordinates - matched case}, there is no mismatch between their histograms and the kernel density estimates. Consequently, the performance of neural network models trained on this data is better, as evidenced by higher precision and recall scores. \\
In this work, we employ a normalization scheme that reduces the mismatch between training and validation data. Intuitively, we think that lesser the mismatch between the distributions of training and validation data, better will be the performance. With such a normalization scheme, the initial conditions for the training data can be sampled from a relatively small interval and the trained neural network models will give reliable performance on different kinds of validation data. \\
As described in section ~\ref{sec:Training and Validation Datasets}, we use $25$ Lorenz systems for training the neural network. We sample the initial conditions of training data according to Eq.~(\ref{subeq : Initial Conditions Train 3}), as we think this interval is a good candidate for training generalizable models. We normalize each of the $25$ Lorenz systems separately. As noted in section ~\ref{sec:Neural Network Architecture}, the input to the neural network is a $6$ dimensional feature vector. We calculate the mean and standard deviation along each of the feature dimensions for each of the Lorenz Systems. Specifically, if $x$ represents x-coordinate of one of the data points of the Lorenz system, we transform $x$ according to Eq.~\ref{subeq:transformation}, where the mean $\mu_{x}$ and standard deviation $\sigma_{x}$ are given by Eq.~(\ref{subeq:mean}) and Eq.~(\ref{subeq:standard deviation}) respectively.
\begin{subequations}
\label{eq:normalization}
    \begin{equation}
    \label{subeq:transformation}
    x = \frac{x - \mu_{x}}{\sigma_{x}}
    \end{equation}
    
    \begin{equation}
    \label{subeq:mean}    
    \mu_{x} = \frac{1}{N}\sum_{i=1}^{N}x_{i}
    \end{equation}
    
    \begin{equation}
    \label{subeq:standard deviation}    
    \sigma_{x} = \frac{1}{N}\sum_{i=1}^{N}(x_{i} - \mu_{x})^2
    \end{equation}
\end{subequations}
We perform normalization for both training and validation data. The Fig.~\ref{Fig : Histogram Comparisons - Before and After Normalization} compares the histograms and kernel density estimates of the training and validation data before and after normalization. One can clearly see that after normalization (Fig.~\ref{fig : Histogram of X-Coordinates - mismatched case after normalization}), the histograms and kernel density estimates of training and validation data have a greater overlap and a reduced mismatch.  

\begin{figure}[!tbp]
    \begin{subfigure}[b]{0.49\columnwidth}
        \includegraphics[width = 0.99\columnwidth]{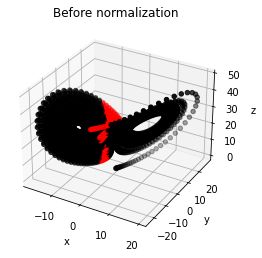}
        \caption{Stable and Unstable Data Points of Lorenz System Before Normalization}
        \label{fig:Stable and Unstable Data Points of Lorenz System Before Normalization}
    \end{subfigure}
    \begin{subfigure}[b]{0.49\columnwidth}
        \includegraphics[width = 0.99\columnwidth]{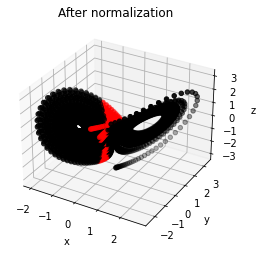}
        \caption{Stable and Unstable Data Points of Lorenz System After Normalization}
        \label{fig : Stable and Unstable Data Points of Lorenz System After Normalization}
    \end{subfigure}
\caption{Comparison of the location of stable and unstable data points of Lorenz System before and normalization. Note that the only change is in the ranges of axes values}
\label{Fig : Lorenz System Comparison before and after normalization}
\end{figure}

The Fig.~\ref{Fig : Lorenz System Comparison before and after normalization} shows the location of stable and unstable data points before and after normalization. Thus, we are able to verify that the normalization scheme does not alter the relative location of stable and unstable data points. It merely shifts the data points along the axes.
\subsection{\label{subsec:Normalization Results}Results}
\begin{figure}[!tbp]
    \begin{subfigure}[b]{0.49\columnwidth}
        \includegraphics[width = 0.99\columnwidth]{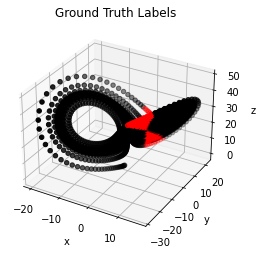}
        \caption{Ground Truth Labels}
        \label{fig:LS_Res6_GT}
    \end{subfigure}
    \begin{subfigure}[b]{0.49\columnwidth}
        \includegraphics[width = 0.99\columnwidth]{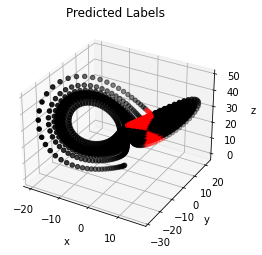}
        \caption{Predicted Labels}
        \label{fig:LS_Res6_Predicted}
    \end{subfigure}
\caption{Classification result using neural network models in mismatched condition after normalization. The initial conditions for training data were sampled according to Eq.~(\ref{subeq : Initial Conditions Train 3}) and for validation data were sampled according to Eq.~(\ref{subeq : Initial Conditions Val 2}). The precision and recall scores are $0.983$ and $0.972$ respectively}
\label{Fig : Res6}
\end{figure}

\begin{figure}[!tbp]
    \begin{subfigure}[b]{0.49\columnwidth}
        \includegraphics[width = 0.99\columnwidth]{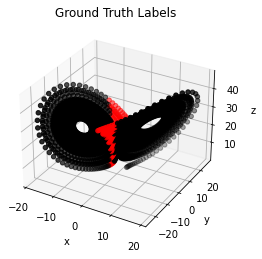}
        \caption{Ground Truth Labels}
        \label{fig:LS_Res7_GT}
    \end{subfigure}
    \begin{subfigure}[b]{0.49\columnwidth}
        \includegraphics[width = 0.99\columnwidth]{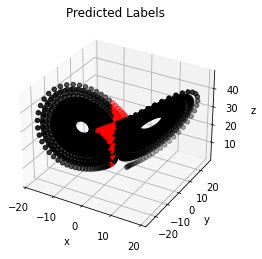}
        \caption{Predicted Labels}
        \label{fig:LS_Res7_Predicted}
    \end{subfigure}
\caption{Classification result using neural network models in mismatched condition after normalization. The initial conditions for training data were sampled according to Eq.~(\ref{subeq : Initial Conditions Train 3}) and for validation data were sampled according to Eq.~(\ref{subeq : Initial Conditions Val 4}). The precision and recall scores are $0.99$ and $0.939$ respectively}
\label{Fig : Res7}
\end{figure}

\begin{figure}[!tbp]
    \begin{subfigure}[b]{0.49\columnwidth}
        \includegraphics[width = 0.99\columnwidth]{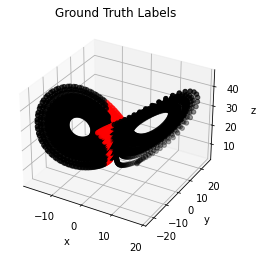}
        \caption{Ground Truth Labels}
        \label{fig:LS_Res8_GT}
    \end{subfigure}
    \begin{subfigure}[b]{0.49\columnwidth}
        \includegraphics[width = 0.99\columnwidth]{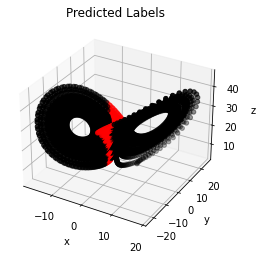}
        \caption{Predicted Labels}
        \label{fig:LS_Res8_Predicted}
    \end{subfigure}
\caption{Classification result using neural network models in mismatched condition after normalization. The initial conditions for training data were sampled according to Eq.~(\ref{subeq : Initial Conditions Train 3}) and for validation data were sampled according to Eq.~(\ref{subeq : Initial Conditions Val 5}). The precision and recall scores are $1$ and $0.91$ respectively}
\label{Fig : Res8}
\end{figure}

\begin{figure}[!tbp]
    \begin{subfigure}[b]{0.49\columnwidth}
        \includegraphics[width = 0.99\columnwidth]{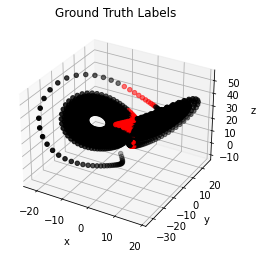}
        \caption{Ground Truth Labels}
        \label{fig:LS_Res9_GT}
    \end{subfigure}
    \begin{subfigure}[b]{0.49\columnwidth}
        \includegraphics[width = 0.99\columnwidth]{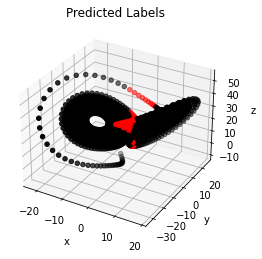}
        \caption{Predicted Labels}
        \label{fig:LS_Res9_Predicted}
    \end{subfigure}
\caption{Classification result using neural network models in mismatched condition after normalization. The initial conditions for training data were sampled according to Eq.~(\ref{subeq : Initial Conditions Train 3}) and for validation data were sampled according to Eq.~(\ref{subeq : Initial Conditions Val 6}). the precision and recall scores are $1$ and $0.95$ respectively.}
\label{Fig : Res9}
\end{figure}

The Fig.~\ref{Fig : Res6} shows the classification result in the mismatched case when the neural network is trained on normalized feature vectors. The initial condition for training are sampled according to Eq.~\ref{subeq : Initial Conditions Train 3} and those for validation are sampled according to Eq. ~\ref{subeq : Initial Conditions Val 2}. The precision are recall scores are $0.983$ and $0.9727$ respectively. One can also see that the results of Fig.~\ref{Fig : Res6} clearly outperform those of Fig.~\ref{fig:Res3}, Fig.~\ref{fig:Res4} and Fig.~\ref{fig:Res5}. This also highlights that neural network models are sensitive to the mismatch between distributions of training and validation data and that normalization schemes are necessary to reduce the mismatch.\\ 
We also verify whether neural networks trained on normalized data , whose initial conditions are sampled according to Eq.~(\ref{subeq : Initial Conditions Train 3}) can perform reliably on a wide variety of validation data. To this end, we test the neural network model on validation data whose initial conditions are sampled according to Eq.~(\ref{subeq : Initial Conditions Val 4}), Eq.~(\ref{subeq : Initial Conditions Val 5}) and Eq.~(\ref{subeq : Initial Conditions Val 6}). The Fig.~\ref{Fig : Res7} shows the result of applying neural network models when the initial conditions of validation data were sampled according to Eq.~(\ref{subeq : Initial Conditions Val 4}). This interval is completely outside of the interval used for sampling the initial conditions of training data (Eq.~(\ref{subeq : Initial Conditions Train 3})). Yet, we see that the classification performance is reasonably good, with a precision score of $0.99$ and a recall score of $0.939$. The Fig.~\ref{Fig : Res8} and Fig. ~\ref{Fig : Res9} show the classification performance when the initial conditions for validation data were sampled according to Eq.~(\ref{subeq : Initial Conditions Val 5}) and Eq.~(\ref{subeq : Initial Conditions Val 6}) respectively. Both these intervals are much wider compared to those for sampling initial conditions for training data (Eq.~(\ref{subeq : Initial Conditions Train 3})). Again, we see that the performance of neural network models is reasonably accurate. The precision and recall scores of the classification result shown in Fig.~\ref{Fig : Res8} are $1$ and $0.91$ respectively. The precision and recall scores of the classification result shown in Fig.~\ref{Fig : Res9} are $1$ and $0.95$ respectively. \\
The Tables \ref{tab:table1} and \ref{tab:table2} summarize the results of our experiments by depicting the average precision and recall scores without and with normalization respectively. For each row, neural network models are trained using training data whose initial conditions are sampled from the intervals specified in the first column. As described in Section \ref{sec:Training and Validation Datasets}, we use 25 Lorenz systems for training and 5 Lorenz systems for validation. The mean precision and recall values are obtained by averaging over the 5 validation lorenz systems. We also show the standard deviation to quantify how much the precision and recall scores deviate from the mean. Per Table \ref{tab:table1}, we observe that the performance of the neural network models drop significantly in the mismatched case (Rows $3$ and $4$ of Table \ref{tab:table1}). We choose the interval $\left[-1,1\right]$ to sample initial conditions for training data for our further experiments, per our observation that it represents a valid subset of the validation data intervals that the neural network is likely to encounter. The results in Table \ref{tab:table2} indicate that our assumption is justified. The neural networks whose training data are sampled from the interval $\left[-1,1\right]$ perform well on a wide variety of validation datasets. These results also show that our normalization scheme greatly helps in improving the performance of neural networks on mismatched data and is a promising step towards training generalizable neural networks for this classification task.
\begin{table*}
\caption{\label{tab:table1}Average Precision and Recall values - Without Normalization}
\begin{ruledtabular}
\begin{tabular}{cccccc}
\multirow{1}{*}{Training} & \multirow{1}{*}{Validation}& \multirow{1}{*}{Mean Precision}&\multirow{1}{*}{Mean Recall}&\multirow{1}{*}{Stddev-Precision}&\multirow{1}{*}{Stddev-Recall}\\ \hline
$\left[0, 1\right]$&$\left[0,1\right]$&$0.82$&$0.617$&$0.122$&$0.170$\\
$\left[-1,0\right]$&$\left[-1,0\right]$&$0.632$&$0.455$&$0.247$&$0.202$\\
$\left[0,1\right]$&$\left[-1,0\right]$&$0.034$&$0.026$&$0.049$&$0.044$\\
$\left[-1,0\right]$&$\left[0,1\right]$&$0.028$&$0.032$&$0.014$&$0.025$\\
$\left[-1,1\right]$&$\left[-1,0\right]$&$0.704$&$0.128$&$0.195$&$0.06$\\
\end{tabular}
\end{ruledtabular}
\end{table*}

\begin{table*}
\caption{\label{tab:table2}Average Precision and Recall values - With Normalization}
\begin{ruledtabular}
\begin{tabular}{cccccc}
\multirow{1}{*}{Training} & \multirow{1}{*}{Interval}& \multirow{3}{*}{Mean Precision}&\multirow{3}{*}{Mean Recall}&\multirow{3}{*}{Stddev-Precision}&\multirow{3}{*}{Stddev-Recall}\\
\multirow{1}{*}{sample} & \multirow{1}{*}{to sample} & {} & {} & {} & {} \\
\multirow{1}{*}{$x_0,y_0,z_0$ from } & \multirow{1}{*}{ics for Validation} &  {} & {} & {} & {} \\
\hline
$\left[-1, 1\right]$&$\left[0,1\right]$&$0.96$&$0.964$&$0.032$&$0.007$\\
$\left[-1,1\right]$&$\left[-1,0\right]$&$0.963$&$0.969$&$0.039$&$0.015$\\
$\left[-1,1\right]$&$\left[-1,1\right]$&$0.982$&$0.978$&$0.007$&$0.005$\\
$\left[-1,1\right]$&$\left[2,4\right]$&$0.988$&$0.934$&$0.01$&$0.015$\\
$\left[-1,1\right]$&$\left[0,10\right]$&$0.975$&$0.955$&$0.025$&$0.028$\\
$\left[-1,1\right]$&$\left[-10,10\right]$&$0.952$&$0.9464$&$0.066$&$0.032$\\
\end{tabular}
\end{ruledtabular}
\end{table*}

\section{\label{sec:Conclusion}Conclusions}
In this paper, we explore the use of neural networks to classify the discrete states of a Lorenz system. Such systems are highly sensitive to the initial conditions and identifying the location of stable and unstable states is a challenging task.Our numeral results suggest that classification performed using a feed forward neural network model can accommodate sensitivity to initial conditions in the training and validation data sets. The classification performance degrades when there is a mismatch between initial conditions used for generating training and validation data sets. We introduce a normalization scheme and show that it significantly improves the classification performance of neural network models  in mismatched conditions. More broadly, our results show the feasibility of using neural networks to study stability aspects of chaotic systems like Lorenz63 systems. Improvement in state estimation via neural network based classification can significantly enhance automated decision making in complex systems within an SDM context. Future work will explore the scalability and explainability aspects of our normalization scheme within neural network based classification applied to large scale complex system state estimation.

\begin{acknowledgments}
The research described in this paper is supported by the Mathematics of Artificial Reasoning in Science (MARS) Initiative at Pacific Northwest National Laboratory (PNNL).  It was conducted under the Laboratory Directed Research and Development Program at PNNL, a multiprogram national laboratory operated by Battelle for the U.S. Department of Energy under  DE-AC05-76RL01830.
\end{acknowledgments}

\bibliography{lorenz}

\end{document}